
\input amstex

\documentstyle{amsppt} 
\magnification1200
\NoBlackBoxes
\pagewidth{6.5 true in}
\pageheight{9.25 true in}
\document

\topmatter 
\title
Weak subconvexity for central values of $L$-functions
\endtitle
\author 
K. Soundararajan
\endauthor
\address
{Department of Mathematics, 450 Serra Mall, Bldg. 380, Stanford University, 
Stanford, CA 94305-2125, USA}
\endaddress
\email 
ksound{\@}stanford.edu
\endemail
\thanks   The author is partially supported by the National Science Foundation (DMS 0500711) 
and the American Institute of Mathematics (AIM). 
\endthanks

\endtopmatter

\def\lam{\lambda}
\def\Lam{\Lambda}
\def\uell{\underline{\ell}}
 \def\uj{\underline{j}}
\document  

\head 1.  Introduction and statement of results \endhead

\noindent A fundamental problem in number theory is to estimate the values of 
$L$-functions at the center of the critical strip.  The Langlands program predicts 
that all $L$-functions arise from automorphic representations of $GL(N)$ over 
a number field, and moreover that such $L$-functions can be decomposed 
as a product of primitive $L$-functions arising from irreducible cuspidal 
representations of $GL(n)$ over ${\Bbb Q}$.   The $L$-functions that 
we consider will either arise in this manner, or will be the Rankin-Selberg 
$L$-function associated to two irreducible cuspidal representations.   Note that 
such Rankin-Selberg $L$-functions are themselves expected to arise from 
 automorphic representations, but this is not known in general.  

Given an irreducible cuspidal automorphic representation $\pi$ (normalized 
to have unitary central character), we  denote the associated 
$L$-function  by $L(s,\pi)$, and its analytic conductor (whose definition we shall 
recall shortly) by $C(\pi)$.  There holds  generally a convexity bound of the form $L(\tfrac 12,\pi) 
\ll_{\epsilon} C(\pi)^{\frac 14+\epsilon}$ (see Molteni [28]).   \footnote{Recently, Roger Heath-Brown [14] has 
pointed out to me an elegant application of Jensen's formula for strips that leads 
generally to the stronger convexity bound $L(\tfrac 12,\pi) \ll C(\pi)^{\frac 14}$.}  
The Riemann hypothesis for $L(s,\pi)$ implies the Lindel{\" o}f hypothesis: $L(\tfrac 12,\pi) \ll C(\pi)^{\epsilon}$.   In several applications it has 
emerged that the convexity bound barely fails to be of use, and that 
any improvement over the convexity bound would have significant 
consequences.   Obtaining such subconvexity bounds has been an 
active area of research, and estimates of the type $L(\tfrac 12,\pi) \ll C(\pi)^{\frac 14-\delta}$ 
for some $\delta >0$ have been obtained for several important classes of $L$-functions.  
However in general the subconvexity problem remains largely open. 
For comprehensive accounts on $L$-functions and the subconvexity problem 
we refer to  Iwaniec and Sarnak [21], and 
Michel [27].

In this paper we describe a method that leads in 
many cases to an improvement over the convexity bound for 
values of $L$-functions.   The improvement is not a saving of a power 
of the analytic conductor, as desired in formulations of 
the subconvexity problem.  Instead we obtain an estimate of the form 
$L(\tfrac 12,\pi) \ll C(\pi)^{\frac 14}/(\log C(\pi))^{1-\epsilon}$, which we term 
{\sl weak} subconvexity. 
In some applications, it appears that a suitable weak subconvexity 
bound would suffice in place of genuine subconvexity.  
In particular, 
by combining sieve estimates 
for the shifted convolution problem developed by Holowinsky together 
with the weak subconvexity 
estimates developed here, Holowinsky and I [18] have been able to resolve 
a conjecture of Rudnick and Sarnak on the mass equidistribution of 
 Hecke eigenforms.


\smallskip
We begin by giving three illustrative examples of our work before 
describing the general result.  

\noindent {\bf Example 1.}   Let $f$ be a holomorphic Hecke eigenform of large weight $k$ 
for the full modular group $SL_2({\Bbb Z})$.  Let $t$ be a fixed real 
number (for example, $t=0$), and consider the symmetric square $L$-function 
$L(\tfrac 12+ it,\text{sym}^2 f)$.   The convexity bound for this 
$L$-function gives $|L(\tfrac 12+it, \text{sym}^2 f)| \ll_t k^{\frac 12+ \epsilon}$, and this 
can be refined to $\ll_t k^{\frac 12}$ by Heath-Brown's remark mentioned in footnote 1.  
Our method gives the weak subconvexity bound, for any $\epsilon >0$,
$$ 
|L(\tfrac 12+it, \text{sym}^2 f)| \ll_{\epsilon} \frac{k^{\frac 12}(1+|t|)^{\frac 34}}{(\log k)^{1-\epsilon}}.  \tag{1.1} 
$$ 
 Obtaining subconvexity bounds in this situation (with a power saving in $k$) remains an 
 important  open problem.  In the case when $k$ is fixed, and $t$ gets large 
 such a subconvexity bound has been achieved recently by  Li [22].  
 We have assumed that the level is $1$ for simplicity, and the result holds for 
 higher level also.  However the assumption that $f$ is holomorphic is 
 essential, since our method makes use of the Ramanujan 
 conjectures known here due to Deligne.  Similar results would hold for Maass forms if 
 we assume the Ramanujan conjectures, but unfortunately the partial results known 
 towards the Ramanujan bounds are insufficient for our purpose.  

\smallskip
\noindent {\bf Example 2.}  Let $f$ be a holomorphic Hecke eigenform of large weight $k$ for 
the full modular group.  Let $\phi$ be a fixed Hecke-Maass eigencuspform for $SL_2({\Bbb Z})$.  
Consider the triple product $L$-function $L(\tfrac 12, f\times f\times \phi)$.  
The convexity bound gives $L(\tfrac 12, f\times f\times \phi) \ll_{\phi} k^{1+\epsilon}$.  
Our weak subconvexity bound gives for any $\epsilon>0$
$$ 
L(\tfrac 12, f\times f\times \phi) \ll_{\phi,\epsilon} \frac{k}{(\log k)^{1-\epsilon}}. \tag{1.2}
$$
Again we could consider higher level, but the assumption that 
$f$ is holomorphic is necessary for our method. 

\smallskip
\noindent {\bf Example 3.}   Let $\pi_0$ be an irreducible cuspidal automorphic representation 
on $GL(m_0)$ over ${\Bbb Q}$ with unitary central character.  We treat $\pi_0$ 
as fixed, and consider $L(\tfrac 12+it, \pi_0)$ in the $t$-aspect.  The convexity 
bound here is $L(\tfrac 12+it,\pi_0) \ll_{\pi_0} (1+|t|)^{\frac{m_0}{4}+\epsilon}$ and 
we obtain 
$$ 
|L(\tfrac 12+it, \pi_0)| \ll_{\pi_0,\epsilon} \frac{(1+|t|)^{\frac{m_0}{4}}}{(\log (1+|t|))^{1-\epsilon}}. \tag{1.3}
$$
Similarly, if $\chi\pmod q$ is a primitive Dirichlet character with $q$ large 
then 
$$ 
L(\tfrac 12,\pi_0 \times \chi) \ll_{\pi_0,\epsilon}  \frac{q^{\frac{m_0}4}}{(\log q)^{1-\epsilon}}. \tag{1.4}
$$
The most general example along these lines is 
the following:  Let $\pi_0$ be as above, and let $\pi$ be an irreducible 
cuspidal automorphic representation on $GL(m)$ with unitary 
central character and {\sl such that $\pi$ satisfies the Ramanujan conjectures}.  
Then we would obtain a weak subconvexity bound for $L(\tfrac 12,\pi_0 \times \pi)$.  

\medskip
We now describe an axiomatic framework (akin to the Selberg class) for 
the class of $L$-functions that we consider.   The properties of $L$-functions 
that we assume are mostly standard, and we have adopted this 
framework in order to clarify the crucial properties needed for our method. 
In addition to the usual assumptions of a Dirichlet series with an Euler product 
and a functional equation, we will need an 
assumption on the size of the Dirichlet series coefficients.  We call this a 
{\sl weak} Ramanujan hypothesis, as the condition is implied by the 
 Ramanujan conjectures.   The reader may prefer to ignore our conditions 
 below and restrict his attention to automorphic $L$-functions satisfying 
 the Ramanujan conjectures, but our framework allows us to deduce 
 results even in cases where the Ramanujan conjectures are not known.  

Let $m\ge 1$ be a fixed natural number.   Let \footnote{Here 
the notation is meant  to suggest that $\pi$ corresponds to an automorphic 
representation, but this is not assumed.} $L(s,\pi)$ be given by the Dirichlet series and 
Euler product
$$ 
L(s,\pi) = \sum_{n=1}^{\infty} \frac{a_{\pi}(n)}{n^s} = \prod_{p}
 \prod_{j=1}^m \Big( 1- \frac{\alpha_{j,\pi}(p)}{p^s} \Big)^{-1}, \tag{1.5a}
$$ 
and we suppose that both the series and product are absolutely convergent in Re$(s)>1$. 
We write 
$$
L(s,\pi_{\infty}) = N^{\frac s2} \prod_{j=1}^{m} \Gamma_{\Bbb R}(s+\mu_j)  \tag{1.5b} 
$$ 
where $\Gamma_{\Bbb R}(s) =\pi^{-s/2}\Gamma(s/2)$, $N$ denotes the conductor, and 
the $\mu_j$ are complex numbers. 
The completed $L$-function $L(s,\pi)L(s,\pi_{\infty})$ has an analytic 
continuation\footnote{Thus we are not allowing $L(s,\pi)$ to have any poles.  It would 
not be difficult to modify our results to allow the completed $L$-function 
to have poles at $0$ and $1$.}  to the entire complex plane, and has finite order.  
Moreover, it  satisfies a functional equation 
$$
 L(s,\pi_{\infty}) L(s,\pi) = \kappa L(1-s,{\tilde \pi}_{\infty} ) L(1-s,\tilde \pi), \tag{1.5c}
$$ 
where $\kappa$ is the root number (a complex number of magnitude $1$), and 
$$ 
L(s,\tilde{\pi}) = \sum_{n=1}^{\infty} \frac{\overline {a_{\pi}(n)}}{n^s}, 
\qquad \text{and } \qquad L(s,\tilde{\pi}_{\infty}) = N^{\frac s2} \prod_{j=1}^{m} \Gamma_{\Bbb R} (s+\overline{\mu_j}).\tag{1.5d}
$$ 
We define the analytic conductor $C=C(\pi)$ (see [21]) by 
$$ 
C(\pi) =N \prod_{j=1}^{m} (1+ |\mu_j|). \tag{1.5e}
$$
Our goal is to obtain an estimate for $L(\tfrac 12,\pi)$ in 
terms of the analytic conductor $C(\pi)$.

Properties (1.5a-d) are standard features of all interesting $L$-functions.  
We now need an assumption on the size of the numbers $\alpha_{j,\pi}(p)$. 
The Ramanujan conjectures, which are expected to hold for 
all $L$-functions, predict that $|\alpha_{j,\pi}(p)| \le 1$ for all $p$.  Further, it 
is expected that the numbers $\mu_j$ appearing in (1.5b) all satisfy Re$(\mu_j)\ge 0$.  
Towards the Ramanujan conjectures it is known (see [26]) that if $\pi$ is an 
irreducible cuspidal representation of $GL(m)$ then $|\alpha_{j,\pi} (p)| \le 
p^{\frac 12- \delta_m}$ for all $p$, and that Re$(\mu_j) \ge -\frac 12+\delta_m$ 
where $\delta_m = 1/(m^2+1)$.  We will make 
the following weak Ramanujan hypothesis.  

Write 
$$ 
-\frac{L^{\prime}}{L}(s,\pi) = \sum_{n=1}^{\infty} \frac{\lam_{\pi}(n)\Lam(n)}{n^s}, \tag{1.6a}
$$ 
where $\lam_\pi(n)=0$ unless $n=p^k$ is a prime power when
it equals $\sum_{j=1}^{m} \alpha_{j,\pi}(p)^k$.  We assume that 
for some constants $A_0, A\ge 1$,  and all $x\ge 1$ there holds 
$$ 
\sum_{x<n \le ex} \frac{|\lam_{\pi}(n)|^2}{n} \Lam(n) \le A^2 + \frac{A_0}{\log ex}.
\tag{1.6b}
$$
Note that the Ramanujan conjecture would give (1.6b) with $A=m$, and 
$A_0 \ll m^2$.  Analogously for the parameters $\mu_j$ we 
assume that\footnote{This assumption is very weak: from [26] we 
know that it holds for all automorphic $L$-functions, and 
also for the Rankin-Selberg $L$-function associated to two automorphic 
representations.}
$$
\text{Re}(\mu_j) \ge -1+\delta_m, \qquad\text{for some } \delta_m >0, \text{ and 
all } 1\le j\le m. \tag{1.6c}
$$


\proclaim{Theorem 1}  Let $L(s,\pi)$ be an $L$-function satisfying 
the properties (1.5a-e) and (1.6a,b,c).  Then for any $\epsilon >0$ 
we have 
$$ 
L(\tfrac 12,\pi ) \ll \frac{C(\pi)^{\frac 14}}{(\log C(\pi))^{1-\epsilon}}.
$$ 
Here the implied constant depends on $m$, $A$, $A_0$, $\delta_m$, and $\epsilon$.  
\endproclaim

We now show how the examples given above fit into the framework of 
Theorem 1. 

\noindent{\bf Example 1 (proof).}   Write the Euler product for $L(s,\text{sym}^2 f)$ as 
$$
L(s,\text{sym}^2 f) = \prod_{p} \Big( 1-\frac{\alpha_p^2}{p^s}\Big)^{-1} 
\Big( 1-\frac{\alpha_p \beta_p}{p^s} \Big)^{-1} \Big( 1-\frac{\beta_p^2}{p^s}\Big)^{-1},
$$
where $\alpha_p = \overline{\beta_p} $ are complex numbers of magnitude $1$ (by Deligne), 
and $\alpha_p+\beta_p$ equals the $p$-th Hecke eigenvalue of $f$.  From the work of 
Shimura we know that the completed $L$-function 
$$
\Lam(s,\text{sym}^2 f) =\Gamma_{\Bbb R}(s+1) \Gamma_{\Bbb R}(s+k-1) \Gamma_{\Bbb R}(s+k) 
L(s,\text{sym}^2 f) 
$$ 
is entire and satisfies the functional equation $\Lam(s,\text{sym}^2 f) =\Lam(1-s,\text{sym}^2 f)$.  
Thus criteria (1.5a-d) hold, and the analytic conductor of (1.5e) is $\asymp k^2$.  
If we write $-L^{\prime}/L(s,\text{sym}^2 f)$ in the notation of (1.6a) then 
the analogous $\lam_{\text{sym}^2 f}(p^k)$ equals $\alpha_p^{2k} +1 + {\beta_p}^{2k }$ 
which is $\le 3$ in magnitude.  Thus criterion (1.6b) holds with $A=3$, and $A_0$ being 
some absolute constant.  Visibly, criterion (1.6c) also holds.  Therefore Theorem 1 applies and yields 
$L(\tfrac 12,\text{sym}^2 f) \ll k^{\frac 12}/(\log k)^{1-\epsilon}$.  

We have shown (1.1) when $t=0$.  To obtain the general case, we apply the 
framework of Theorem 1 to the shifted function $L_t(s,\text{sym}^2 f):= L(s+it,\text{sym}^2 f)$.  
The criteria we require hold.   The only difference is that we must make corresponding 
shifts to the $\Gamma$-functions appearing in the functional equation.  These 
shifts imply that the analytic conductor is now $\asymp (1+|t|)(k+|t|)^2 \ll k^2 (1+|t|)^3$.  
Applying Theorem 1, we complete the proof of (1.1). 

\smallskip 

\noindent{\bf Example 2 (proof).}     Write the $p$-th Hecke eigenvalue of 
$f$ as $\alpha_{f}(p) + \beta_f(p)$ where $\alpha_f(p)\beta_f(p)=1$ 
and $|\alpha_f(p)|=|\beta_f(p)| =1$.  Write the $p$-th Hecke eigenvalue 
of $\phi$ as $\alpha_{\phi}(p) +\beta_\phi (p)$ where $\alpha_{\phi}(p)\beta_{\phi}(p)=1$, 
but we do not know here the Ramanujan conjecture that these are both of size $1$.  Write 
also the Laplace eigenvalue of $\phi$ as $\lam_\phi = \frac 14+ t_{\phi}^2$, 
where\footnote{This is true since we are working on the full modular group.  
For a congruence subgroup, we could use Selberg's bound that the 
least eigenvalue is $\ge \frac 3{16}$ which 
gives that $|\text{Im}(t_{\phi})|\le \frac 14$; see [26].} $t_{\phi}\in {\Bbb R}$.  

The triple product $L$-function $L(s,f\times f\times \phi)$ is then defined by 
means of the Euler product of degree $8$ (absolutely convergent in Re$(s)>1$) 
$$ 
\align
\prod_p &\Big( 1- \frac{\alpha_{f}(p)^2 \alpha_{\phi}(p)}{p^s}\Big)^{-1} 
\Big( 1- \frac{\alpha_{\phi}(p)}{p^s}\Big)^{-2} \Big(1-\frac{\beta_f(p)^2 \alpha_{\phi}(p)}{p^s}\Big)^{-1} 
\\
&\times \Big(1-\frac{\alpha_f(p)^2 \beta_\phi(p)}{p^s}\Big)^{-1}\Big( 1-\frac{\beta_\phi(p)}{p^s}\Big)^{-2} 
\Big(1-\frac{\beta_f(p)^2 \beta_\phi(p)}{p^s}\Big)^{-1}.
\\
\endalign
$$
This $L$-function is not primitive and factors as $L(s,\phi) L(s,\text{sym}^2 f\times \phi)$.  
Consider the product of eight $\Gamma$-factors 
$$ 
L_{\infty}(s,f\times f\times \phi):= \prod_{\pm} \Gamma_{\Bbb R}(s+k-1\pm it_{\phi}) 
\Gamma_{\Bbb R}(s+k \pm it_{\phi}) \Gamma_{\Bbb R}(s\pm it_{\phi}) 
\Gamma_{\Bbb R}(s+1\pm it_{\phi}).
$$
From the work of Garrett [6], it is known that $L(s,f\times f \times \phi) L_{\infty}(s,f\times f\times \phi)$ 
is an entire function in ${\Bbb C}$, and its value at $s$ equals its value at $1-s$.  
Thus criteria (1.5a-d) are met, and the analytic conductor in (1.5e) is of 
size $\ll k^4 (1+|t_{\phi}|)^8$.  

If we write $-\frac{L^{\prime}}{L}(s,f) =\sum_{n} \lam_f(n) \Lam(n)n^{-s}$ and $-\frac{L^{\prime}}L(s,\phi)
= \sum_{n} \lam_\phi(n) \Lam(n) n^{-s}$, then 
$$
-\frac{L^{\prime}}{L}(s,f\times f\times \phi) = \sum_{n} \frac{\lam_{f\times f\times \phi}(n)\Lam(n)}{n^s} 
= \sum_{n} \frac{\lam_f(n)^2 \lam_\phi(n)\Lam(n)}{n^s}.
$$
If $n=p^k$ then $|\lam_f(n)| = |\alpha_f(p)^k +\beta_f(p)^k| \le 2$, and so to check (1.6b) 
for $f\times f\times \phi$, we need only show that 
$$ 
\sum_{x < n\le ex} \frac{|\lam_{\phi}(n)|^2 \Lam(n)}{n} \le A^2 +\frac{A_0}{\log (ex)},
$$ 
for all $x\ge 1$, where $A$ and $A_0$ are constants which {\sl are allowed to 
depend on  $\phi$}.   This condition follows from an appeal to the 
Rankin-Selberg theory for $L(s,\phi \times \phi)$ which is known to 
extend analytically to ${\Bbb C}$ except for a simple pole at $s=1$.  Since 
$-\frac{L^{\prime}}{L}(s,\phi\times \phi) = \sum_{n} |\lam_\phi(n)|^2 \Lam(n) n^{-s}$, 
and $L(s,\phi \times \phi)$ has a classical zero-free region $\text{Re}(s) > 1-c_{\phi}/\log (1+|t|)$ 
(see Theorem 5.44 of [20]), 
we may deduce, arguing as in the proof of the prime number theorem, that 
$$ 
\sum_{ x< n\le ex} \frac{|\lam_{\phi}(n)|^2 \Lam(n)}{n} = 1 +O_\phi\Big(\frac{1}{\log (ex)}\Big),
$$ 
from which our desired weak Ramanujan estimate follows.   Criterion (1.6c) is immediate 
from our formula for the $\Gamma$ factors above.

Thus Theorem 1 applies and we obtain the desired estimate (1.2).   We could 
also obtain weak subconvexity bounds for triple products $\phi \times f_1 \times f_2$, 
fixing $\phi$ (a holomorphic or Maass eigencuspform) 
and allowing $f_1$ and $f_2$ to vary over holomorphic eigencuspforms.  
Or we could consider $f_1\times f_2 \times f_3$ with all three varying over 
holomorphic eigencuspforms.  Subconvexity bounds (with a power saving) have been 
obtained for triple products $\phi_1 \times \phi_2\times \phi_3$ where 
$\phi_1$ and $\phi_2$ are considered fixed, and $\phi_3$ varies over 
holomorphic or Maass forms, see [1] and [32].  
 
\smallskip 

\noindent{\bf Example 3 (proof).}  This example follows upon using the 
ideas in the proof of Examples 1 and 2.  The weak Ramanujan 
hypothesis (1.6b) is verified by an appeal to the Rankin-Selberg theory 
for $L(s,\pi_0\times {\tilde \pi_0})$ as in Example 2.   The general Rankin-Selberg 
theory is the culmination of work by many authors, notably Jacquet, Piatetskii-Shapiro, 
and Shalika, Shahidi, and Moeglin and Waldspurger; a convenient synopsis 
of the analytic features of this theory may be found in [30].    A narrow 
zero-free region (which is sufficient to check (1.6b)) for general Rankin-Selberg $L$-functions 
has been established by Brumley [3]. 
The value $L(\tfrac 12+it, \pi_0)$ is bounded by shifting $L$-functions 
as in Example 1.  We omit further details. 

\smallskip

\noindent{\bf Application to the mass equidistribution of Hecke eigenforms.} 
Perhaps the most interesting application of our results pertains to a conjecture of 
Rudnick and Sarnak on the mass equidistribution of Hecke eigenforms.  For simplicity, 
consider a holomorphic Hecke eigencuspform $f$ of weight $k$ for the full 
modular group $\Gamma = SL_2({\Bbb Z})$.   Consider the measure 
$$ 
\mu_f = y^k |f(z)|^2 \frac{dx \ dy}{y^2},
$$ 
where we suppose that $f$ has been normalized to 
satisfy 
$$ 
\int_{\Gamma \backslash {\Bbb H}} y^k |f(z)|^2 \frac{dx \ dy}{y^2} = 1.
$$
Rudnick and Sarnak ([29], see also [25, 31]) have 
conjectured that as $k\to \infty$, the measure $\mu_f$ 
approaches the uniform distribution measure $\frac 3\pi \frac{dx dy}{y^2}$ on the 
fundamental domain $X=\Gamma\backslash {\Bbb H}$.   This is the 
holomorphic analog of their quantum unique ergodicity conjecture for Maass forms.  Lindenstrauss [23]
has made great progress on the latter question, but his ergodic theoretic 
methods do not seem to apply to the holomorphic case.  

Set $F_k(z) = y^{k/2} f(z)$, and recall the Petersson inner product of two nice functions 
$g_1$ and $g_2$ on $X$ 
$$ 
\langle g_1, g_2 \rangle  = \int_{X} g_1(z) \overline{g_2(z)} \frac{dx \ dy}{y^2}. 
$$
A smooth bounded function on $X$ has a spectral 
expansion in terms of the constant function, the space of Maass cusp forms which 
are eigenfunctions of the Laplacian and all Hecke operators, and the 
Eisenstein series $E(z,\tfrac 12+it)$ with $t\in {\Bbb R}$, see Iwaniec [19].  
  Thus by an analog of Weyl's equidistribution criterion, the Rudnick-Sarnak conjecture 
amounts to showing that 
$$ 
\langle\phi F_k, F_k  \rangle,\qquad \langle E(\cdot, \tfrac 12+it)F_k, F_k \rangle 
 \to 0
$$ 
as $k\to \infty$, where $\phi$ is a fixed Maass cusp form which is an eigenfunction of 
the Laplacian and all Hecke operators, and $E(z,\tfrac 12 +it)$ denotes the Eisenstein 
series (with $t$ fixed). 

Using the unfolding method,  it is easy to show that 
$$
|\langle  E(\cdot ,\tfrac 12+it) F_k, F_k\rangle | 
= \Big|\pi^{\frac 32}
 \frac{\zeta(\tfrac 12+it) L(\tfrac 12+it, \text{sym}^2 f)} {\zeta(1+2it) L(1,\text{sym}^2 f)} 
 \frac{\Gamma(k-\tfrac 12+it)}{\Gamma(k)} \Big|.
 $$
Since $|\Gamma(k-\tfrac 12+it)| \le \Gamma(k-\tfrac 12)$, $|\zeta(\tfrac 12+it)| \ll (1+|t|)^{\frac 14}$, 
and $|\zeta(1+2it)| \gg 1/\log (1+|t|)$,  using Stirling's formula and 
our bound (1.1) it follows that  
 $$ 
|\langle   E(\cdot,\tfrac 12+it)F_k, F_k\rangle|  \ll_{\epsilon} \frac{(1+|t|)^{2}}{(\log k)^{1-\epsilon} L(1,\text{sym}^2 f)}. 
 $$
 
  For the case of a Maass cusp form, a beautiful formula of Watson (see Theorem 3 of [33]) shows that (here $\phi$ has been normalized so that $\langle \phi, \phi \rangle =1$)
  $$ 
  |\langle \phi F_k, F_k \rangle|^2 
  = \frac 18 \frac{L_{\infty}(\tfrac 12, f\times f\times \phi) L(\tfrac 12, f\times f\times \phi)}{\Lam(1,\text{sym}^2 f)^2 \Lam(1, 
  \text{sym}^2 \phi)  } 
  $$
  where  
  $L(s,f\times f\times \phi)$ is the triple product $L$-function of Example 2, 
  and $L_{\infty}$ denotes its Gamma factors (see Example 2 (proof)), and 
  $$
  \Lam(s, \text{sym}^2 f) = \Gamma_{\Bbb R}(s+1)\Gamma_{\Bbb R}(s+k-1) \Gamma_{\Bbb R}(s+k) L(s,
  \text{sym}^2 f), 
  $$
  and 
  $$ 
  \Lam(s,\text{sym}^2 \phi) = \Gamma_{\Bbb R}(s) \Gamma_{\Bbb R} (s+ 2it_{\phi}) \Gamma_{\Bbb R} 
  (s-2it_{\phi}) L(s, \text{sym}^2 \phi). 
  $$
  Using Stirling's formula and the bound (1.2) of Example 2 we conclude that 
  $$ 
  |\langle \phi F_k, F_k \rangle | \ll_{\phi} \frac{1}{(\log k)^{\frac 12-\epsilon} L(1,\text{sym}^2 f)}. 
  $$
  
\proclaim{Corollary 1} With notations as above, we have 
$$ 
|\langle \phi F_k, F_k  \rangle | \ll_{\phi, \epsilon} \frac{1}{(\log k)^{\frac 12-\epsilon}L(1,\text{\rm sym}^2 f)}. 
$$   
Moreover 
$$ 
|\langle E(\cdot, \tfrac 12+it) F_k, F_k \rangle | \ll_{\epsilon} \frac{(1+|t|)^2 }{(\log k)^{1-\epsilon}L(1,
\text{\rm sym}^2 f)}.
$$ 
\endproclaim

Given $\delta >0$, Corollary 1 shows that if $f$ 
ranges over those Hecke eigencuspforms with 
$L(1,\text{\rm sym}^2 f) \ge (\log k)^{-\frac 12+\delta}$ 
then as $k\to \infty$, the measure $\mu_f$ converges to $\frac 3\pi \frac{dx \ dy}{y^2}$.
This criterion on $L(1,\text{sym}^2 f)$ is expected to hold for all eigenforms $f$; 
for example, it is implied by the Riemann hypothesis 
for $L(s,\text{sym}^2 f)$. 
Using large sieve estimates\footnote{Precisely, with very few exceptions 
one can approximate $L(1,\text{sym}^2 f)$ by a short Euler product.  
Such results in the context of Dirichlet $L$-functions are classical, 
and a detailed account may be found in [7].  In the context of symmetric 
square $L$-functions, see Luo [24], and Cogdell and Michel [4].}  one can show 
that the number of exceptional eigenforms $f$ with weight $k\le K$ for 
which the criterion fails is $\ll K^{\epsilon}$.  
Our criterion complements the work of Holowinsky [17], who attacks the 
mass equidistribution conjecture by an entirely different method.  Combining 
his results with ours gives a complete resolution of the Rudnick-Sarnak 
conjecture on mass equidistribution for eigenforms.  A detailed account of 
this result will appear 
in a joint work with Holowinsky [18].

   As noted earlier, 
the barrier to using our methods for Maass forms is the Ramanujan conjecture 
which remains open here.  The weak Ramanujan hypothesis that 
we need could be verified (using the large sieve) for all but $T^{\epsilon}$ 
Maass forms with Laplace eigenvalue below $T$.  Thus with at most $T^{\epsilon}$ 
exceptions, one could establish the equidistribution of Maass forms on $SL_2({\Bbb Z}) 
\backslash {\Bbb H}$.  

For simplicity we have confined ourselves to $\Gamma=SL_2({\Bbb Z})$ above. 
Similar results apply to congruence subgroups of level $N$.
One could also 
consider the cocompact case of quaternion division algebras.  If the quaternion algebra 
is unramified at infinity, Lindenstrauss's work shows the equidistribution of Maass 
cusp forms of large eigenvalue.  Our results would show the corresponding 
equidistribution for holomorphic eigenforms, allowing for a small 
number of exceptional cases (at most $K^{\epsilon}$ exceptions with 
weight below $K$).  In the case of a ramified quaternion algebra (acting on the 
unit sphere $S^2$), the problem concerns the equidistribution of eigenfunctions on the 
sphere (see [2]), and again our results establish such equidistribution 
except for a small number of cases (omitting at most $N^{\epsilon}$ 
spherical harmonics of degree below $N$ that are also Hecke eigenforms).  In these 
compact cases, there is no analog of Holowinsky's work, and so we are unable 
to obtain a definitive result as in $SL_2({\Bbb Z})\backslash {\Bbb H}$.

\medskip







We now return to the setting of Theorem 1, and describe an 
auxiliary result which will be used to prove Theorem 1.  
For an $L$-function satisfying (1.5a-e) and (1.6a,b,c), we may use  
the convexity bound to establish that (see Lemma 4.2 below)\footnote{We recall here 
that $L(s,\pi)$ was assumed not to have any poles.  If we alter our framework to allow 
a pole at $s=1$, say, then (1.7) would be modified to an asymptotic formula 
with a main term of size $x$.  Then Theorem 2 would extrapolate that asymptotic 
formula to a wider region.} 
$$ 
 \sum_{n\le x} a_{\pi}(n) \ll  \frac{x}{\log x}, \tag{1.7} 
$$ 
provided $x\ge C^{\frac 12} (\log C)^{B}$ for some positive constant $B$.  
Our main idea is to show 
that similar cancellation holds even when $x= C^{\frac 12} (\log C)^{-B}$ 
for any constant $B$.  
 
 \proclaim{Theorem 2}  Let $L(s,\pi)$ be as in Theorem 1.   For any $\epsilon>0$, any 
 positive constant $B$, and all $x\ge C^{\frac 12}(\log C)^{-B}$ we have
 $$ 
 \sum_{n\le x} a_{\pi} (n) \ll \frac{x}{(\log x)^{1-\epsilon}}.
 $$ 
 The implied constant may depend on $A$, $A_0$, $m$, $\delta_m$, $B$ and $\epsilon$. 
 \endproclaim

Once Theorem 2 is established, Theorem 1 will follow from a standard 
partial summation argument using an approximate functional equation. 
In Theorems 1 and 2, by keeping track of the various parameters involved, 
it would be possible to quantify $\epsilon$.   However, the limit of 
our method would be to obtain a bound $C^{\frac 14}/\log C$ in 
Theorem 1, and $x/\log x$ in Theorem 2.  

We have termed the saving of $(\log C)^{1-\epsilon}$ as weak 
subconvexity, and as noted above this is close to the 
limit of our method.  One may legitimately call a saving of 
any fixed power of $\log C$ as weak subconvexity.  For example, 
in the application to the Rudnick-Sarnak conjecture any log power 
saving together with Holowinsky's work would suffice to show 
that $|\langle\phi F_k ,F_k \rangle | \to 0$, for $\phi$ a fixed Maass cusp form.   
However to deal with the Eisenstein series contributions in that 
application, a saving of a substantial power of $\log C$ is needed\footnote{A saving of $(\log C)^{0.7}$ 
would probably be sufficient.}, and 
a small saving would not suffice.   
Finally, it would be very desirable to establish a version of weak subconvexity 
saving a large power of $\log C$.  For example if one could save $(\log C)^{2+\delta}$ for 
any fixed $\delta >0$,  
one would obtain immediately the mass equidistribution consequences 
mentioned above.  Similarly it would be desirable to improve upon the 
weak Ramanujan criterion that we have imposed.



  {\bf Acknowledgments.}  I am grateful to Peter Sarnak for encouragement 
  and several valuable discussions.   The ideas in this paper build upon my 
  work with Andrew Granville on multiplicative functions; it is a pleasure 
  to thank him for the many stimulating discussions we have had on that 
  topic.  I would also like to thank Philippe Michel and Akshay Venkatesh 
  who referred me to Holowinsky's work, which led to the joint work [18].  Thanks also 
  to Tom Watson for clarifying his formula, and to Roman Holowinsky 
  for a careful reading of this paper.

\head 2.  Slow oscillation of mean values of multiplicative functions \endhead

\noindent We now discuss the ideas underlying Theorem 2, and state our 
main technical result from which the results stated in \S 1 follow.   At the 
heart of Theorem 2 is the fact that mean values of multiplicative functions 
vary slowly.  Knowing (1.7) in the range $x\ge C^{\frac 12}(\log C)^{B}$, 
this fact will enable us to extrapolate (1.7) to the range 
$x\ge C^{\frac 12}(\log C)^{-B}$.  

The possibility of obtaining such extrapolations was first considered 
by Hildebrand [15, 16].  If $f$ is a multiplicative function, we shall denote by 
$S(x)=S(x;f)$ the partial sum $\sum_{n\le x}f(n)$.  
Hildebrand [16] showed that if $-1\le f(n) \le 1$ is a real 
valued multiplicative function then for $1\le w\le \sqrt{x}$
$$ 
\frac{1}{x} \sum_{n\le x} f(n)  = \frac{w}{x} \sum_{n\le x/w} f(n) + O \Big( \Big( \log \frac{\log x}{\log 2w}\Big)^{-\frac 12}\Big). \tag{2.1}
$$
In other words, the mean value of $f$ at $x$ does not change 
very much from the mean-value at $x/w$.  Hildebrand [15] used this 
idea to show that from knowing Burgess's character sum estimates\footnote{For simplicity, 
suppose that $q$ is cube-free.}
for $x\ge q^{\frac 14+\epsilon}$ one may obtain some 
non-trivial cancellation even in the range $x\ge q^{\frac 14-\epsilon}$.  

Elliott [5] generalized Hildebrand's work to cover complex valued 
multiplicative functions with $|f(n)| \le 1$, and also strengthened 
the error term in (2.1).   Notice that a direct extension of (2.1) for 
complex valued functions is false.  Consider $f(n)=n^{i\tau}$ 
for some real number $\tau \neq 0$.  Then $S(x;f) 
= x^{1+i\tau}/(1+i\tau) + O(1)$, and $S(x/w;f) = (x/w)^{1+i\tau}/(1+i\tau)+ O(1)$.  Therefore (2.1) 
is false, and instead we have that $S(x)/x$ is close to $w^{i\tau} S(x/w)/(x/w)$.  
Building on the pioneering work of Halasz [11, 12] on mean-values of multiplicative functions, 
Elliott showed that for a multiplicative function $f$ with $|f(n)| \le 1$, 
there exists a real number $\tau = \tau(x)$ with $|\tau| \le \log x$ 
such that for $1\le w\le \sqrt{x}$ 
$$ 
S(x) = w^{1+i\tau} S(x/w) + O\Big( x\Big(\frac{\log 2w}{\log x}\Big)^{\frac 1{19}}\Big). \tag{2.2} 
$$ 
In [8], Granville and Soundararajan give variants and stronger versions 
of (2.2), with $\tfrac 1{19}$ replaced by $1-2/\pi -\epsilon$. 

In order to establish Theorem 2, we require similar results when the 
multiplicative function is no longer constrained to the unit disc.  The situation 
here is considerably more complicated, and instead of showing that 
a suitable linear combination of $S(x)/x$ and $S(x/w)/(x/w)$ is small, 
we will need to consider linear combinations involving 
several terms $S(x/w^j)/(x/w^j)$ with $j=0$, $\ldots$, $J$.  In order to 
motivate our main result, it is helpful to consider two illustrative 
examples.  
 
\noindent {\bf Example 2.1.} Let $k$ be a natural number, and take $f(n) = d_k(n)$, the $k$-th
divisor function.  Then, it is easy to show that  $S(x) = x P_k(\log x) + O(x^{1-1/k+\epsilon})$ 
where $P_k$ is a polynomial of degree $k-1$.  If $k\ge 2$, it follows that 
$S(x)/x -S(x/w)/(x/w)$ is of size $(\log w)(\log x)^{k-2}$, which is not $o(1)$.  
However, if $1\le w\le x^{1/2k}$,  the linear combination 
$$ 
\sum_{j=0}^{k} (-1)^j \binom{k}{j} \frac{S(x/w^j)}{x/w^j} =  \sum_{j=0}^{k} (-1)^j \binom{k}{j} P_k(\log x/w^j) 
+ O(x^{-\frac{1}{2k}}) = O(x^{-\frac{1}{2k}})
$$
is very small.  

\noindent {\bf Example 2.2.}  Let $\tau_1$, $\ldots$, $\tau_R$ be distinct real numbers, and 
let $k_1$, $\ldots$, $k_R$ be natural numbers.   Let $f$ be the multiplicative function 
defined by $F(s) =\sum_{n=1}^{\infty} f(n)n^{-s} = \prod_{j=1}^{R} \zeta(s-i\tau_j)^{k_j}$.  
Consider here the linear combination (for $1\le w\le x^{1/(2(k_1+\ldots+k_R))}$) 
$$ 
\frac 1 x 
\sum_{j_1=0}^{k_1} \cdots \sum_{j_R=0}^{k_R} (-1)^{j_1+ \ldots+j_R} \binom{k_1}{j_1}\cdots
\binom{k_R}{j_R} w^{j_1(1+i\tau_1)+\ldots+ j_R(1+i\tau_R)} S\Big(\frac{x}{w^{j_1+\ldots+j_R}}\Big). 
$$
By Perron's formula we may express this as, for $c>1$, 
$$ 
\frac{1}{2\pi i} \int_{c-i\infty}^{c+i\infty} \prod_{j=1}^{R} \zeta(s-i\tau_j)^{k_j} 
(1-w^{1+i\tau_j-s})^{k_j} x^{s-1} \frac{ds}{s}.
$$ 
Notice that the poles of the zeta-functions at $1+i\tau_j$ have been cancelled by the factors $(1-w^{1+i\tau_j-s})^{k_j}$.
Thus the integrand has a pole only at $s=0$, and a standard contour shift argument 
shows that this integral is $\ll x^{-\delta}$ for some $\delta >0$.  

Fortunately, it turns out that Example 2.2 captures the behavior of mean-values of the multiplicative functions of interest to us.  In order to state our result, we require some notation.  Let 
$f$ denote a multiplicative function and recall that 
$$ 
S(x) =S(x;f) = \sum_{n\le x} f(n). \tag{2.3}
$$ 
We shall write 
$$ 
F(s) = \sum_{n=1}^{\infty} \frac{f(n)}{n^s}, \tag{2.4}
$$ 
and we shall assume that this series converges absolutely in Re$(s)>1$.  
Moreover we write 
$$ 
-\frac{F^{\prime}}{F}(s) = \sum_{n=1}^{\infty} \frac{\lam_f(n)\Lam(n)}{n^s} 
= \sum_{n=1}^{\infty}\frac{ \Lam_f(n)}{n^s}, \tag{2.5}
$$ 
where $\lam_f(n)=\Lam_f(n)=0$ unless $n$ is the power of a prime $p$.  We next assume 
the analog of the weak Ramanujan hypothesis (1.6b).  Namely, we 
suppose that there exist constants $A$, $A_0 \ge 1$ such that for all $x\ge 1$ 
we have 
$$ 
\sum_{x< n\le ex} \frac{|\lam_f(n)|^2 \Lam(n)}{n} \le A^2 +\frac{A_0}{\log (ex)}. 
\tag{2.6}
$$

Let $R$ be a natural number, and let $\tau_1$, $\ldots$, $\tau_R$ denote 
$R$ real numbers.  Let ${\uell} = (\ell_1,\ldots,\ell_R)$ and ${\uj}=(j_1,\ldots,j_R)$ 
denote vectors of non-negative integers, with the notation ${\uj} \le {\uell}$ 
indicating that $0\le j_1\le \ell_1$, $\ldots$, $0\le j_R\le \ell_R$.  
Define 
$$ 
\binom{\uell}{\uj} = \binom{\ell_1}{j_1} \cdots \binom{\ell_R}{j_R}. \tag{2.7}
$$
Finally, we define a measure of the oscillation of the mean-values of $f$ by 
setting 
$$ 
\align
{\Cal O}_{\uell}(x,w) &= {\Cal O}_{\uell}(x,w;\tau_1,\ldots,\tau_R) \\
&= \sum_{\uj\le \uell} (-1)^{j_1+\ldots+j_R} \binom{\uell}{\uj}  
w^{j_1(1+i\tau_1)+\ldots+j_R(1+i\tau_R)}S\Big( \frac{x}{w^{j_1+\ldots+j_R}}\Big).
\tag{2.8}
\endalign
$$

\def\uL{\underline{L}}
\proclaim{Theorem 2.1} Keep in mind the conditions and notations (2.3) through (2.8).  
Let $X\ge 10$ and $1\ge \epsilon >0$ be given. Let $R= [10A^2/\epsilon^2]+1$ and 
put $L=[10AR]$, and $\uL = (L,\ldots,L)$.   Let $w$ be such that $0\le \log w \le 
(\log X)^{\frac{1}{3R}}$.  There exist real numbers $\tau_1$, $\ldots$, 
$\tau_R$ with $|\tau_j| \le \exp((\log \log X)^2)$ such that 
for all $2\le x\le X$ we have 
$$ 
|{\Cal O}_{\uL}(x,w;\tau_1,\ldots,\tau_R)| \ll \frac{x}{\log x} (\log X)^{\epsilon}.
$$
The implied constant above depends on $A$, $A_0$ and $\epsilon$. 
\endproclaim

For a general multiplicative function, we cannot hope for any better 
bound for the oscillation than $x/\log x$.  To see this, suppose $w\ge 2$, 
and consider the multiplicative function $f$ with $f(n)=0$ for 
$n\le x/2$ and $f(p)=1$ for primes $x/2<p\le x$.  Then $S(x) \gg x/\log x$ 
whereas $S(x/w^j) =1$ for all $j\ge 1$, and therefore for any choice 
of the numbers $\tau_1$, $\ldots$, $\tau_R$ we would 
have ${\Cal O}_{\uL}(x,w) \gg x/\log x$.  

Our proof of Theorem 2.1 builds both on the techniques of 
Halasz (as developed in [5] and [8]), and also the idea of {\sl pretentious} multiplicative 
functions developed by Granville and Soundararajan (see [9] and [10]).   In \S 5
we will describe the choice of the numbers $\tau_1$, $\ldots $, $\tau_R$ appearing 
in Theorem 2.1, and develop relevant estimates for the Dirichlet series $F(s)$.  
Then the proof of Theorem 2.1 will be completed in \S 6.

\head 3. Some preliminary Lemmas \endhead 

\noindent We collect together here some Lemmas that 
will be useful below.  We begin with a combinatorial lemma. 

\proclaim{Lemma 3.1}  Let $b(1)$, $b(2)$, $\ldots$ be a sequence of complex numbers.  
Define the sequence $a(0)=1$, $a(1)$, $a(2)$, $\ldots$ by means 
of the formal identity 
$$ 
\exp\Big( \sum_{k=1}^{\infty} \frac{b(k)}{k} x^k\Big) = \sum_{n=0}^{\infty} a(n) x^n.  \tag{3.1}
$$ 
For $j=1$ or $2$, define the sequences $A_j(0)=1$, $A_j(1)$, $A_j(2)$, $\ldots$ by means of the 
formal identity 
$$ 
\exp\Big( \sum_{k=1}^{\infty} \frac{|b(k)|^j}{k} x^k \Big) = \sum_{n=0}^{\infty} A_j(n) x^n. 
$$ 
Then $A_j(n) \ge |a(n)|^j$ for all $n$. 
\endproclaim
\demo{Proof}  If we expand out the LHS of (3.1) and equate coefficients we obtain that 
$$
a(n) = \sum\Sb \underline{\lam} \\ \lam_1+\ldots+\lam_r=n \\ \lam_1 \ge \lam_2 \ge \ldots \ge \lam_r 
\ge 1\endSb b(\lam_1)\cdots b(\lam_r) W(\underline{\lam}), \tag{3.2}
$$
where the sum is over all partitions ${\underline{\lam}}$ of $n$, and $W({\underline{\lam}})>0$ 
is a weight attached to each partition which is independent of the sequence $b(k)$.  Although we can write down explicitly what the 
weight $W(\underline{\lam})$ is, we do not require this.  All we need is that 
$$ 
\sum_{\underline{\lam}} W(\underline{\lam}) = 1, 
$$ 
which follows from (3.2) by taking $b(k)=1$ for all $k$ so that $\exp(\sum_{k=1}^{\infty} b(k)x^k/k) = 
1/(1-x) =\sum_{n=0}^{\infty} x^n$ whence $a(n)=1$.   When $j=1$ the Lemma follows by the 
triangle inequality.  In the case $j=2$, Cauchy-Schwarz gives 
$$ 
\align
|a(n)|^2 &\le \Big( \sum_{\underline{\lam}} |b(\lam_1) \cdots b(\lam_{r})|^2 W(\underline{\lam}) \Big) 
\Big( \sum_{\underline {\lam}} W({\underline{\lam}}) \Big) \\
&=  \Big( \sum_{\underline{\lam}} |b(\lam_1) \cdots b(\lam_{r})|^2 W(\underline{\lam}) \Big)= A_2(n).\\
\endalign
$$ 
\enddemo

The significance of Lemma 3.1 for us is the following.  Let $f$ be a multiplicative 
function with the Euler factor at $p$ being (compare (1.6a) and (2.5)) 
$$ 
\sum_{n=0}^{\infty}\frac{ f(p^n)}{ p^{ns}} = \exp\Big( \sum_{k=1}^{\infty} \frac{\lam_f(p) }{kp^{ks}} \Big)  
$$
then the Lemma guarantees that if we write 
$$ 
\exp\Big( \sum_{k=1}^{\infty} \frac{|\lam_f(p)|^2}{k p^{ks}} \Big) = \sum_{n=0}^{\infty} \frac{f^{(2)}(p^n)}{p^{ns}}, 
$$ 
and use this to define a multiplicative function $f^{(2)}(n)$, then $|f(n)|^2 \le f^{(2)}(n)$.  
In the case that $f(n)=a_{\pi}(n)$ corresponds to the coefficients of an automorphic $L$-function $\pi$, 
this means that the coefficient\footnote{We assume that $n$ is the power of an unramified prime.} 
 $a_{\pi\times {\tilde\pi}}(n)$ of the Rankin-Selberg $L$-function exceeds $|a_{\pi}(n)|^2$.   The reader may compare this Lemma with a 
similar result of Molteni, Proposition 6 of  [28].  Although we are only interested in Lemma 3.1 in 
the cases $j=1$ and $2$, an application of H{\" o}lder's inequality shows that 
a similar result holds for all $j\ge 1$.

\proclaim{Lemma 3.2}  Let $f$ be a multiplicative function, and keep the notations 
(2.4) and (2.5), and suppose that the criterion (2.6) holds.  For all $x\ge 2$ we have 
$$ 
\sum_{n\le x} \frac{|f(n)|}{n} \ll (\log x)^{A}, \qquad \text{and} \qquad 
\sum_{n\le x} \frac{|f(n)|^2}{n} \ll (\log x)^{A^2}.
$$
Moreover, for all $2\ge \sigma>1$ we have 
$$
|F(\sigma+it)| \ll \Big(\frac{1}{\sigma-1}\Big)^A \qquad \text{and} \qquad 
|F^{\prime}(\sigma+it)| \ll \Big(\frac{1}{\sigma-1}\Big)^{A+1}. 
$$ 
The implied constants may depend on $A$ and $A_0$. 
\endproclaim
\demo{Proof}  These are simple consequences of our weak Ramanujan 
  assumption (2.6).  By splitting the sum into intervals $e^{k} < n \le e^{k+1}$ and using 
 (2.6), we see that for any $2 \ge \sigma >1$ 
 $$ 
 \sum_{n=2}^{\infty} \frac{|\lam_f(n)|^2 \Lam(n)}{n^{\sigma} \log n} 
 \le A^2 \log\Big( \frac{1}{\sigma-1}\Big) + O(1), \tag{3.3} 
 $$ 
 where the error term above depends on $A$ and $A_0$.  
  By Cauchy-Schwarz it follows also 
  that 
  $$ 
  \sum_{n=2}^{\infty} \frac{|\lam_f(n)|\Lam(n)}{n^{\sigma} \log n} \le A \log \Big(\frac{1}{\sigma-1}\Big) + O(1).
 \tag{3.4}
 $$
 
 Using Lemma 3.1 and (3.4) we see that 
 $$ 
 \sum_{n\le x} \frac{|f(n)|}{n} \ll \sum_{n=1}^{\infty} \frac{|f(n)|}{n^{1+1/\log (ex)} }
 \le \exp\Big( \sum_{n\ge 2} \frac{|\lam_f(n)|\Lam(n)}{n^{1+1/\log (ex)}\log n}\Big) 
 \ll (\log x)^A.
 $$
 This proves our first inequality.  The second inequality follows in the 
 same way, using (3.3) instead of (3.4).  The third inequality follows easily from (3.4). 
  Finally, for $2\ge \sigma >1$,
  $$
  \Big| \frac{F^{\prime}}{F}(\sigma+it)\Big| \le \sum_{n\ge 2} 
  \frac{|\lam_f(n)|\Lam(n)}{n^{\sigma} } \ll \frac{1}{\sigma-1} 
  $$
  using Cauchy-Schwarz and (2.6) as in the proof of (3.4).  The Lemma follows.

\enddemo 

\proclaim{Lemma 3.3} Let $f$ be a multiplicative function as in Lemma 3.2.  
Then for all $x\ge 2$ 
$$ 
\sum_{n\le x} |f(n)| \ll x(\log x)^A.$$
Moreover, for $1\le y\le x$ we have 
$$ 
\sum_{x<n\le x+y}| f(n)| \ll (yx)^{\frac 12} (\log x)^{A^2/2}.
$$
\endproclaim

\demo{Proof}  Since $\sum_{n\le x} |f(n)| \le x\sum_{n\le x} |f(n)|/n$, 
the first assertion follows from Lemma 3.2.  Cauchy-Schwarz gives 
$$ 
\Big|\sum_{x<n\le x+y} |f(n)| \Big|^2\le y \sum_{x<n\le x+y} |f(n)|^2 
\ll yx \sum_{n\le 2x} \frac{|f(n)|^2}{n},
$$ 
and our second assertion also follows from Lemma 3.2.

\enddemo

 \head 4.  Deduction of the main results from Theorem 2.1 \endhead 
 
 \noindent In this section we shall show how Theorems 1 and 2 
 follow from Theorem 2.1.   We begin with Theorem 2, whose proof 
 will require the following simple convexity bound for our $L$-functions.  
 
 \proclaim{Lemma 4.1} Let $L(s,\pi)$ be an $L$-function satisfying 
 the properties (1.5a-e) and (1.6a,b,c).  Then for all $t \in {\Bbb R}$ we 
 have 
 $$ 
 |L(\tfrac 12+it, \pi)| \ll C(\pi)^{\frac 14} (1+|t|)^{\frac m4 +1} (\log C(\pi))^{A}.
 $$
 \endproclaim 
 
 \demo{Proof}   Define $\Lam(s)=L(s,\pi) L(s,\pi_{\infty}) e^{(s-\frac 12-it)^2}$. 
Using the Phragmen-Lindel{\" of} principle 
 we may bound $|\Lam(\tfrac 12+it)|$ by the maximum value taken by $|\Lam(s)|$ 
 on the lines Re$(s)=1+1/\log C$, and Re$(s)=-1/\log C$.   The functional equation 
 shows that the maximum on the line Re$(s)=-1/\log C$ is the same as the maximum on 
 the line Re$(s)=1+1/\log C$.  Therefore, using Lemma 3.2,
 $$ 
 \align
 |L(\tfrac 12+it,\pi)| &\le  \max_{y\in {\Bbb R}} \frac{|\Lam(1+1/\log C+it+iy)|}{|L(\tfrac 12+it,\pi_{\infty})|} 
\\
& \ll (\log C)^{A} \max_{y\in {\Bbb R}} e^{-y^2} \frac{|L(1+1/\log C+it+iy,\pi_{\infty})|}{|L(\tfrac 12+it,\pi_{\infty})|}.\\
\endalign
 $$
Using Stirling's formula, we may show that 
$$ 
\Big| \frac{\Gamma_{\Bbb R}(1+1/\log C+it+iy+\mu_j)}{\Gamma_{\Bbb R}(\tfrac 12+it+\mu_j)}\Big| 
\ll e^{2|y|}  (1+|t|+|\mu_j|)^{\frac 14+\frac 1{2\log C}},
$$
 where we used that Re$(\mu_j) \ge -1+\delta_m$ to ensure that the numerator stays away from poles 
 of the $\Gamma$-function.   The Lemma follows immediately.
 \enddemo
 
 Our next Lemma establishes the result (1.7) stated in the Introduction. 
 
 \proclaim{Lemma 4.2}  Let $L(s,\pi)$ be as above.  In the range 
 $x\ge x_0:= C(\pi)^{\frac 12} (\log C(\pi))^{50 mA^2}$ we have 
 $$ 
 \sum_{n\le x} a_{\pi}(n) \ll \frac{x}{\log x}.
 $$
 \endproclaim 
 
 \demo{Proof}  Observe that for any $c>0$, $y>0$, any $\lam >0$, and any 
natural number $K$,    
$$ 
\align
\frac{1}{2\pi i} \int_{c-i\infty}^{c+i\infty} 
\frac{y^s}{s} \Big( \frac{e^{\lam s}-1}{\lam s} \Big)^K ds 
&= \frac{1}{\lam^K} \int_0^{\lam} \ldots \int_0^{\lam} \frac{1}{2\pi i}\int_{c-i\infty}^{c+i\infty} (ye^{x_1+\ldots +x_K})^s \frac{ds}{s} 
dx_1 \ldots dx_K \\ 
&= \cases 
1 &\text{if } y \ge 1 \\ 
\in [0,1] &\text{if } 1 >y \ge e^{-\lam K}\\
0 &\text{if } y<e^{-\lam K}. \\
\endcases
\\
\endalign
$$ 
Therefore, for any $c>1$, 
$$ 
\frac{1}{2\pi i} \int_{c-i\infty}^{c+i\infty} L(s,\pi) \frac{x^s}{s} \Big( \frac{e^{\lam s}-1}{\lam s}\Big)^{K} 
ds = \sum_{n\le x} a_{\pi}(n) + O\Big( \sum_{x<n\le e^{K\lam } x } |a_{\pi}(n)|\Big). \tag{4.1} 
$$
We shall take here $K=[m/4]+3$ and $\lam = (\log x)^{-2-A^2}$.  
Then for large $x$ we have $e^{K\lam} \le 2$, and so applying Lemma 3.3 to the multiplicative function $a_{\pi}(n)$ we see that  
the error term above is 
$$ 
\align
&\ll (e^{K\lam}-1)^{\frac 12} x (\log x)^{A^2/2} \ll \frac{x}{\log x}. \tag{4.2}
\\
\endalign
$$ 

Now we move the line of integration in the LHS of (4.1) to the line $c=1/2$.   
Using Lemma 4.1 we see that 
the integral on the $1/2$ line is 
$$ 
\ll C(\pi)^{\frac 14} x^{\frac 12} \lam^{-K} (\log C)^{A} 
\int_{-\infty}^{\infty} (1+|t|)^{m/4+1} \frac{dt}{(1+|t|)^{K+1}} 
\ll C(\pi)^{\frac 14} x^{\frac 12} (\log (Cx))^{12 m A^2}.
$$ 
 Combining this 
 with (4.1) and (4.2) we conclude that for 
 $x\ge x_0:= C(\pi)^{\frac 12} (\log C(\pi))^{50 mA^2}$ we have 
 $$ 
 \sum_{n\le x} a_{\pi}(n) \ll \frac{x}{\log x}, 
  $$
 proving our Lemma.

 \enddemo

\demo{Proof of Theorem 2}  To prove Theorem 2, we invoke Theorem 2.1.  Let 
$R=[10A^2/\epsilon^2]+1$ 
and $L=[10AR]$ be as in Theorem 2.1.  Let $x_0$ be as in Lemma 4.2, and 
let $x_0 \ge x\ge C^{\frac 12}/(\log C)^B$.  Take $w= x_0/x$ and $X=xw^{LR}$.  Applying Theorem 2.1 to 
the multiplicative function $a_\pi$ (note that (1.6b) gives the assumption (2.6)) we 
find that for an appropriate choice of $\tau_1$, $\ldots$, $\tau_R$ that
$$ 
|{\Cal O}_{\uL} (X,w) | \ll \frac{X}{(\log X)^{1-\epsilon}}. \tag{4.3} 
$$
But, by definition, the LHS above is  
$$
w^{LR} \Big|\sum_{n\le X/w^{LR}} a_\pi(n) \Big|+ 
O\Big( \sum_{j=0}^{LR-1} w^j \Big| \sum_{n\le X/w^j} a_{\pi}(n)\Big|\Big). \tag{4.4} 
$$
Now $X/w^{LR} = x$, and for $0\le j\le LR-1$ we have 
$X/w^j \ge xw = x_0$ so that the bound of Lemma 4.2 applies.  Therefore (4.4) equals 
$$ 
 w^{LR} \Big|\sum_{n\le x} a_{\pi}(n)\Big| 
+ O\Big( \frac{X}{\log X}\Big),
$$
From (4.3) we conclude that 
$$ 
\Big| \sum_{n\le x} a_{\pi}(n)\Big| \ll w^{-LR} \frac{X}{(\log X)^{1-\epsilon}}  
\ll \frac{x}{(\log x)^{1-\epsilon}},
$$
which proves Theorem 2.

 \enddemo
 
\demo{Deduction of Theorem 1 from Theorem 2}   Theorem 1 follows from 
Theorem 2 by a standard argument using an ``approximate functional 
equation" for $L(\tfrac 12,\pi)$ (see for example Harcos [13], Theorem 2.1) 
and partial summation.  For the sake of completeness we provide a 
brief argument.  We start with, for $c>\frac 12$, 
$$ 
\frac{1}{2\pi i} \int_{c-i\infty}^{c+i\infty} L(s+\tfrac 12,\pi)\frac{L(s+\frac 12,\pi_{\infty})}{L(\frac 12,\pi_{\infty})} 
e^{s^2}\frac{ds}{s},
$$ 
and move the line of integration to Re$(s)=-c$.  We encounter a pole at $s=0$, and 
so the above equals 
$$ 
L(\tfrac 12,\pi) + \frac{1}{2\pi i} \int_{-c-i\infty}^{-c+i\infty} L(s+\tfrac 12,\pi) \frac{L(s+\frac 12,\pi_{\infty})}{L(\frac 12,\pi_{\infty})} 
e^{s^2}\frac{ds}{s}.
$$ 
Now we use the functional equation above, and make a change of variables $s\to -s$.  
In this way we obtain that 
$$ 
\align
L(\tfrac 12,\pi) &= \frac{1}{2\pi i} \int_{c-i\infty}^{c+i\infty} L(s+\tfrac 12,\pi)\frac{L(s+\frac 12,\pi_{\infty})}{L(\frac 12,\pi_{\infty})} 
e^{s^2}\frac{ds}{s} \\
&\hskip .5 in + \frac{\kappa}{2\pi i} \int_{c-i\infty}^{c+i\infty} L(s+\tfrac 12,{\tilde \pi})\frac{L(s+\frac 12,{\tilde \pi}_{\infty})}{L(\frac 12,{ \pi}_{\infty})} e^{s^2} \frac{ds}{s}.
\\
\endalign
$$

Consider the first integral above; the second is estimated similarly.  Using 
$$ 
L(s+\tfrac 12,\pi) = (s+\tfrac 12) \int_1^{\infty} \sum_{n\le x} a_{\pi}(n) \frac{dx}{x^{s+\frac 32}},
$$ 
we see that the first integral above equals 
$$ 
\int_1^{\infty} \sum_{n\le x} a_{\pi}(n) \Big( \frac{1}{2\pi i} \int_{c-i\infty}^{c+i\infty} (s+\tfrac 12) \frac{L(s+\frac12, \pi_\infty)}{L(\frac 12,\pi_{\infty})} e^{s^2} x^{-s}\frac{ds}{s}
\Big) \  \frac{dx}{x^{\frac 32}}. \tag{4.5}
$$
To estimate the inner integral over $s$, we move the line of integration 
either to Re$(s)=\frac 12-\frac{\delta_m}{2}$, or to Re$(s)=2$.  Using 
Stirling's formula, we see that this inner integral is $\ll \min ( (\sqrt{C}/x)^{\frac 12-\frac{\delta_m}{2}}, (\sqrt{C}/x)^2)$.  Thus (4.5) is 
$$ 
\ll C^{\frac 14-\frac{\delta_m}{4}} \int_1^{\sqrt{C}} \Big| \sum_{n\le x} a_{\pi}(n) \Big| 
\frac{dx}{x^{2-\frac{\delta_m}{2}} }+ C \int_{\sqrt{C}}^{\infty} \Big| \sum_{n\le x} a_{\pi}(n)\Big| 
\frac{dx}{x^{\frac 72}}.
$$
We now split into the ranges $x\le \sqrt{C}/(\log C)^{4A/\delta_m}$ and 
$x>\sqrt{C}/(\log C)^{4A/\delta_m}$. 
In the first range we use Lemma 3.3 to bound $|\sum_{n\le x} a_{\pi}(n)|$ by $\ll x(\log x)^A$, 
and in the second range we use $\sum_{n\le x} a_{\pi}(n) \ll x/(\log x)^{1-\epsilon}$ by Theorem 2. 
Inserting these bounds above, we conclude that the quantity in (4.5) is $\ll C^{\frac 14}/(\log C)^{1-\epsilon}$, and Theorem 1 follows.

\enddemo 
 
 
 
 
 
 
\head 5.  Successive maxima \endhead 

\noindent  Recall the conditions and notations (2.3) through (2.8).   
As in Theorem 2.1, $X\ge 10$ and $1\ge \epsilon >0$ are given, and 
 $R= [10A^2/\epsilon^2]+1$.  In this section we define the points $\tau_1$, 
 $\ldots$, $\tau_R$ appearing in Theorem 2.1, and 
 collect together some estimates for the Dirichlet series $F(s)$.  
 
 From now on, we shall write $T= \exp((\log \log X)^2)$.  We define 
 $\tau_1$ to be that point $t$ in the compact set ${\Cal C}_1 = [-T,T]$ where the maximum of 
 $|F(1+1/\log X+it)|$ is attained.   Now remove the interval 
 $(\tau_1-(\log X)^{-\frac 1R},\tau_1+(\log X)^{-\frac 1R})$ 
 from ${\Cal C}_1= [-T,T]$, and let ${\Cal C}_2$ denote the 
 remaining compact set.   We define $\tau_2$ to be that
 point $t$ in ${\Cal C}_2$ where the 
 maximum of $|F(1+1/\log X+it)|$ is attained.  
 Next remove the interval $(\tau_2-(\log X)^{-\frac 1R},\tau_2+(\log X)^{-\frac 1R})$ 
 from ${\Cal C}_2$ leaving behind the compact set ${\Cal C}_3$.  Define 
 $\tau_3$ to be the point where the maximum of $|F(1+1/\log X+it)|$ for
 $t\in {\Cal C}_3$ is attained.  We proceed in this manner, defining the 
 successive maxima $\tau_1$, $\ldots$, $\tau_R$, and the nested compact sets ${\Cal C}_1 \supset 
 {\Cal C}_2 \supset \ldots \supset {\Cal C}_R$.  Notice that 
 all the points $\tau_1$, $\ldots$, $\tau_R$ lie in $[-T,T]$, and moreover are well-spaced: 
 $|\tau_j-\tau_k| \ge (\log X)^{-\frac 1R}$ for $j\neq k$.  
 
  

Lemma  3.2 bounds $|F(1+1/\log X+it)|$ by  $\ll (\log X)^A$.  For $t\in [-T,T]$ 
we will show that a much better bound holds, unless $t$ happens to be 
near one of the points $\tau_1$, $\ldots$, $\tau_R$.  The next Lemma 
is inspired by the ideas in [9] and [10].

\proclaim{Lemma 5.1}  Let $1\le j\le R$ and let $t$ be a point in ${\Cal C}_j$.  
Then 
$$ 
|F(1+1/\log X+it)| \ll (\log X)^{A \sqrt{1/j+(j-1)/(jR)}}.
$$
In particular if $t \in {\Cal C}_R$ 
we have $|F(1+1/\log X+it)| \ll (\log X)^{\epsilon/2}$.
\endproclaim
\demo{Proof} If $t \in {\Cal C}_j$ then for all $1\le r\le j$ 
$$
|F(1+1/\log X+it)| \le |F(1+1/\log X+i\tau_j)| \le |F(1+1/\log X+i\tau_r)|.
$$ 
Therefore, 
$$ 
\align
|F(1+1/\log X+i\tau_j)| &\le \Big( \prod_{r=1}^j |F(1+1/\log X+i\tau_r)|\Big)^{\frac 1j} \\
&\le \exp\Big( \text{Re} \frac{1}{j} \sum_{n\ge 2} \frac{\lam_{f}(n)\Lam(n)}{n^{1+1/\log X}(\log n)} (n^{-i\tau_1}+\ldots +n^{-i\tau_j})
\Big).\\
\endalign
$$
By Cauchy-Schwarz 
$$ 
\align
\sum_{n\ge 2} \frac{|\lam_{f}(n)|\Lam(n)}{n^{1+1/\log X}\log n} \Big| \sum_{r=1}^{j} n^{-i\tau_r} \Big| 
&\le \Big( \sum_{n\ge 2} \frac{|\lam_{f}(n)|^2\Lam(n)}{n^{1+1/\log X}\log n} \Big)^{\frac 12} 
\\
&\hskip 1 in \times \Big( \sum_{n\ge 2} \frac{\Lam(n)}{n^{1+1/\log X}\log n} \Big| \sum_{r=1}^{j} 
n^{-i\tau_r} \Big|^2\Big)^{\frac 12}.\\
\endalign
$$ 
By (2.6) the first factor above is $\le (A^2 \log \log X + O(1))^{\frac 12}$.  
To handle the second factor, we expand out the square and obtain 
$$ 
\align 
\sum_{n\ge 2} \frac{\Lam(n)}{n^{1+1/\log X}\log n} &\Big| \sum_{r=1}^{j} n^{-i\tau_r}\Big|^2 \\
&= j \sum_{n\ge 2} \frac{\Lam(n)}{n^{1+1/\log X}\log n} + 2 \text{Re }\sum\Sb 1\le r <s \le j\endSb \sum_{n\ge 2} 
\frac{\Lam(n)}{n^{1+1/\log X +i(\tau_r-\tau_s)}\log n} \\
&= j (\log \log X+ O(1)) +2 \sum_{1\le r< s\le j} \log |\zeta(1+1/\log X+i(\tau_r-\tau_s))|.\\ 
\endalign
$$
Now note that $(\log X)^{-\frac{1}{R}} \le |\tau_r-\tau_s| \le 2T$ and hence 
$|\zeta(1+1/\log X+i(\tau_r-\tau_s))| \le (\log X)^{\frac 1R} +O(1)$.   Using this above, 
the Lemma follows.

\enddemo

Let $\uell= (\ell_1,\ldots,\ell_R)$ be a vector of non-negative integers.  
In our proof of Theorem 2.1 we will encounter (recall Example 2.2 from \S 2 
where a similar quantity arises)
$$ 
{\Cal F}_{\uell}(s) = F(s) \prod_{j=1}^{R} (1- w^{1+i\tau_j -s} )^{\ell_j}. \tag{5.1}
$$ 
We will need good bounds for this quantity, and we record such 
estimates in the next two Lemmas.  

\proclaim{Lemma 5.2}  Let $\sigma \ge 1+1/\log X$.  Then 
$$ 
\max_{|t|\le T/2} |{\Cal F}_{\uell} (\sigma+it)| \le \max_{|t|\le T} |{\Cal F}_{\uell}(1+1/\log X+it)| 
+ O((\log X)^{-1}).
$$
\endproclaim 
\demo{Proof} This follows from the argument of Lemma 2.2 in Granville and Soundararajan [8].  
For completeness we give a proof.  Put $\sigma=1+1/\log X+\alpha$, and assume that 
$\alpha>0$.   The Fourier transform of $k(z)=e^{-\alpha|z|}$ is ${\hat k}(\xi) = 
\int_{-\infty}^{\infty} e^{-\alpha |z|-i\xi z} dz = \frac{2\alpha}{\alpha^2+\xi^2}$.  By 
Fourier inversion, we have for any $z\ge 1$
$$
\align
z^{-\alpha} &= k(\log z) =k(-\log z)= \frac{1}{2\pi }\int_{-\infty}^{\infty} {\hat k}(\xi) z^{-i\xi} d\xi \\
&=\frac{1}{\pi } \int_{-T/2}^{T/2} \frac{\alpha}{\alpha^2+\xi^2} z^{-i\xi}d\xi + O\Big( \frac{\alpha}{T}\Big). 
\\
\endalign
$$ 
Using this relation appropriately we 
obtain that 
$$ 
{\Cal F}_{\uell}(\sigma+it) = \frac{1}{\pi} \int_{-T/2}^{T/2} \frac{\alpha}{\alpha^2+\xi^2} 
{\Cal F}_{\uell}(1+1/\log X+it+i\xi) d\xi +O\Big(\frac{\alpha}{T} \sum_{n=1}^{\infty} \frac{|f(n)|}{n^{1+1/\log X}} \Big).
$$
Using Lemma 3.2 and partial summation, the error term above is $O((\log X)^A/T) = O((\log X)^{-1})$. 
If $|t|\le T/2$ then $|t+\xi|\le T$, and so 
$$
\max_{|t|\le T/2} |{\Cal F}_{\uell}(\sigma+it)| \le \max_{|y|\le T} |{\Cal F}_{\uell}(1+1/\log X+iy)| 
\frac{1}{\pi} \int_{-T/2}^{T/2} \frac{\alpha}{\alpha^2+\xi^2} d\xi + O((\log X)^{-1}),
$$ 
and the Lemma follows.
\enddemo 


\proclaim{Lemma 5.3}  Suppose $\ell_j \ge L-1$ for all $1\le j\le R$ where we recall 
that $L=[10AR]$.  Then provided 
$0\le \log w \le (\log X)^{1/(3R)}$ we have 
$$ 
\max_{|t|\le T} |{\Cal F}_{\uell}(1+1/\log X+it)| \ll (\log X)^{\epsilon/2}.
$$
\endproclaim
 \demo{Proof} Suppose first that $|t|\le T$ but $|t-\tau_j|> (\log X)^{-1/R}$ 
 for all $1\le j \le R$.  Then Lemma 5.1 gives that $|F(1+1/\log X+ it)| \ll (\log X)^{\epsilon/2}$ and 
 so $|{\Cal F}_{\uell}(1+1/\log X+ it)| \ll (\log X)^{\epsilon/2}$ as well. 
 
 Now suppose that $|t-\tau_j| \le (\log X)^{-1/R}$ for some $1\le j\le R$.  By Lemma 3.2 we have 
 that 
 $|F(1+1/\log X+it)|\ll (\log X)^{A}$.
 Moreover  
 $$
 |1-w^{-1/\log X-it+i\tau_j}|^{\ell_j} \ll \Big(\frac{\log w}{(\log X)^{\frac 1R}} \Big)^{\ell_j} 
 \le (\log X)^{-\frac{2(L-1)}{3R}}\ll (\log X)^{-A},
 $$
 and hence  $|{\Cal F}_{\uell}(1+1/\log X+it)| \ll 1$.
  The Lemma follows.
 
 \enddemo

Using our work so far, we can record a preliminary estimate for the oscillation which 
we shall refine in the next section to obtain Theorem 2.1. 

\proclaim{Proposition 5.4}  Suppose $\ell_j \ge L-1$ for all $1\le j\le R$, and that $0\le \log w \le 
(\log X)^{\frac{1}{3R}}$.  For $x\le X$ we have 
$$ 
{\Cal O}_{\uell}(x,w) \ll x (\log X)^{2\epsilon/3}.
$$
\endproclaim 
\demo{Proof}  Since $S(x) \ll x(\log x)^A$ by Lemma 3.3, we may assume that 
$\log x \ge (\log X)^{\epsilon/(2A)}$, and in particular $x\ge w^{2RL}$ is large.  
 By Perron's formula we have that for $c>1$ 
$$ 
\frac{1}{2\pi i} \int_{c-i\infty}^{c+i\infty} F(s) z^s \Big(\frac{e^{s/\sqrt{T}}-1}{s/\sqrt{T}} \Big) \frac{ds}{s} 
= \sum_{n\le z} f(n) + O\Big( \sum_{z  < n\le ze^{1/\sqrt{T}} }|f(n)|\Big).
$$  
By Lemma 3.3, the error term above is easily seen to be $O(z/\log z)$ in the 
range $T \le z\le X$.  
Using the above formula in the definition of the oscillation, we obtain that 
for $w^{2RL} \le x\le X$
$$ 
{\Cal O}_{\uell}(x,w) = 
\frac{1}{2\pi i} \int_{c-i\infty}^{c+i\infty} F(s) \prod_{j=1}^{R} (1-w^{1+i\tau_j-s})^{\ell_j} 
\Big(\frac{e^{s/\sqrt{T}}-1}{s/\sqrt{T}} \Big) \frac{ds}{s} +O\Big( \frac{x}{\log x}\Big). 
$$
We choose $c=1+1/\log X$ and split the integral into two parts: when $|{\text{Im}}(s)| \le T$ 
and when $|{\text {Im}}(s)|>T$.  For the first range we use Lemma 5.3, and 
so this portion of the integral contributes 
$$ 
\align
&\ll x  (\log X)^{\epsilon/2}   \int_{|\text{Im}(s)|\le T} \frac{|ds|}{|s|}
\ll x  (\log X)^{2\epsilon/3}.
\\
\endalign
$$
For the second region we use that $|F(s)|\ll (\log X)^A$ and deduce that this 
integral contributes 
$$ 
\ll x (\log X)^A \int_{|\text{Im}(s)|>T} \frac{\sqrt{T}}{|s|} \frac{|ds|}{|s|} \ll x.
$$
The Proposition follows.
\enddemo

\head 6. Proof of Theorem 2.1 \endhead

\def\Lam{\Lambda}

\noindent In this section we shall prove Theorem 2.1, with the points 
$\tau_1$, $\ldots$, $\tau_R$ being the successive maxima defined in \S 5.  Recall the conditions and 
notations (2.3) through (2.8), and the notation introduced in \S 5.  
Recall that $R= [10A^2/\epsilon^2]+1$, that $L=[10AR]$ and $\uL=(L,\ldots,L)$.   
Throughout we assume that $\log w \le (\log X)^{\frac 1{3R}}$, that $x\le X$,  
and all implicit constants will be allowed to depend on $A$, $A_0$ and $\epsilon$. 

\def\ue{\underline{e}}

\proclaim{Lemma 6.1}  
With the above notations, we have 
$$ 
\align
(\log x) {\Cal O}_{\uL}(x,w) &= 
\sum_{d\le x} \Lam_f(d) {\Cal O}_{\uL}(x/d,w) + O(x(\log X)^{\epsilon}).\\
\endalign
$$
\endproclaim

\demo{Proof}  Write $\log x =\log (x/w^{j_1+\ldots+j_R}) + (j_1+\ldots+j_R) \log w$.  Hence, we may 
express $(\log x ){\Cal O}_{\uL}(x,w)$ as 
$$ 
\align
&\sum\Sb\uj \le \uL \endSb  (-1)^{j_1+\ldots+j_R} 
\binom{\uL}{\uj} \log (x/w^{j_1+\ldots+j_R}) S(x/w^{j_1+\ldots+j_R}) w^{j_1(1+i\tau_1)+\ldots+j_R (1+i\tau_R)} 
\\
&  + \log w \sum_{\uj \le \uL} (-1)^{j_1+\ldots+j_R} (j_1+\ldots+j_R)  \binom{\uL}{\uj} 
S(x/w^{j_1+\ldots+j_R}) w^{j_1(1+i\tau_1)+\ldots+j_R(1+i\tau_R)}. \\
\endalign
$$
The second term above is readily seen to be 
$$ 
- \sum_{k=1}^{R} L w^{1+i\tau_k} \log w \ {\Cal O}_{\uL-\ue_k}(x/w,w ), \tag{6.1}
$$
where we let ${\ue_k}$ denote the vector with $1$ in the $k$-th place and 
$0$ elsewhere.  Note that the coordinates of $\uL-\ue_k$ are all at least $L-1$, and 
so by Proposition 5.4, the quantity in (6.1) is $\ll x(\log X)^{\epsilon}$.  

To analyze the first term, we write 
$$ 
\align
(\log x)  S(x) &= \sum_{n\le x} f(n)  \log n + \sum_{n\le x} f(n) \log (x/n)\\ 
&=\sum_{n\le x} \sum_{d|n} \Lam_f(d) f(n/d) + \int_1^x S(t) \frac{dt}{t}\\
&= \sum_{d\le x} \Lam_f(d) S(x/d) + \int_1^x S(t) \frac{dt}{t}.\\
\endalign
$$ 
Therefore the first term equals 
$$ 
\sum_{d\le x} \Lam_f(d) {\Cal O}_{\uL}(x/d,w ) + \int_1^x {\Cal O}_{\uL} (t,w)\frac{dt}{t}, 
$$ 
where we used that $\int_1^{x/w^j} S(t) dt/t = \int_{1}^{x} S(t/w^j) dt/t$. 
By Proposition 5.4, the integral above is $\ll x (\log X)^{\epsilon}$.  The Lemma 
follows.
 
 \enddemo
 

 \proclaim{Lemma 6.2}  
For $1\le z\le y$ with $y+z\le X$ we have  
$$ 
\align
\Big| |{\Cal O}_{\uL}(y,w)|^2 - |{\Cal O}_{\uL}(y+z,w)|^2 \Big|&
 \ll y(\log X)^{\epsilon} \sum_{j=0}^{LR} w^j  \sum_{y/w^j < n\le (y+z)/w^j} |f(n)|.\\
\endalign
$$
 \endproclaim 
 \demo{Proof}  The quantity we wish to estimate is 
 $$ 
  \le 
 \Big( |{\Cal O}_{\uL}(y,w)| + |{\Cal O}_{\uL}(y+z,w)|\Big) \Big| {\Cal O}_{\uL}(y+z,w) - {\Cal O}_{\uL} 
 (y,w)\Big|.
 $$  
 By Proposition 5.4, the first factor is $\ll y (\log X)^{\epsilon}$.  The second factor above 
 is 
 $$ 
 \ll \sum_{j=0}^{LR}   w^j \Big| S((y+z)/w^j) -S(y/w^j)\Big| \ll 
 \sum_{j=0}^{LR} w^j \sum_{y/w^j < n\le (y+z)/w^j} |f(n)|. 
 $$ 
 
 \enddemo

\proclaim{Proposition 6.3}  
We have
$$ 
\align
\log x \  |{\Cal O}_{\uL}(x,w)| &\ll x(\log \log x)^{\frac 12} \Big( \int_1^x \log (ey) |{\Cal O}_{\uL}(y,w)|^2 
\frac{dy}{y^3} \Big)^{\frac 12} + x(\log X)^{\epsilon}.\\
\endalign
$$
\endproclaim 
\demo{Proof}   We start with Lemma 6.1,  
and are faced with estimating $\sum_{d\le x} |\Lam_f(d)| |{\Cal O}_{\uL}(x/d,w)|$. 
We split this sum into the terms $d\le D:= [\exp((\log \log X)^6)]$ and $d>D$.  
For the first category of terms we use Proposition 5.4 and obtain that this contribution is 
$$ 
\sum_{d\le D} |\Lam_f(d)| \frac{x}{d} (\log X)^{2\epsilon/3} \ll x (\log X)^{\epsilon}, 
$$ 
upon using (2.6). 

It remains to estimate the contribution of the terms $d>D$.  Define temporarily the 
function $g(t) = t \log (ex/t)$ for $1\le t\le x$.  
By Cauchy-Schwarz we have 
$$
\align
 \sum_{D<d\le x} |\Lam_f(d)| |{\Cal O}_{\uL}(x/d,w)| &\le \Big( \sum_{D<d \le x} 
 \frac{| \lam_f(d)|^2 \Lam(d)}{g(d)} \Big)^{\frac 12} 
\Big ( \sum_{D<d\le x} g(d) \Lam (d) |{\Cal O}_{\uL}(x/d,w)|^2\Big)^{\frac 12}\\
&\ll (\log \log x)^{\frac 12}\Big(\sum_{D<d\le x} g(d)\Lam(d) |{\Cal O}_{\uL}(x/d,w)|^2 \Big)^{\frac 12},
\tag{6.2}\\
\endalign
$$
 where the last estimate follows from (2.6) and partial summation.

Put $\psi_0(x) =\sum_{n\le x} (\Lam(n)-1) 
 = \psi(x)-x$ so that $\psi_0(x) =O(x\exp(-c\sqrt{\log x}))$ by the prime number theorem.  
 Then 
 $$ 
  \align 
  \sum_{D< d\le x} g(d) \Lam (d) |{\Cal O}_{\uL}(x/d,w)|^2& =   \sum_{D<d \le x} g(d)  |{\Cal O}_{\uL}(x/d,w)|^2 \\
&  + 
  \sum_{D< d \le x} 
(\psi_0(d)-\psi_0(d-1))g(d) |{\Cal O}_{\uL}(x/d,w)|^2. 
\tag{6.3} \\
\endalign
$$

We may rewrite the second term in the RHS of (6.3) as 
$$ 
\align
\sum_{D<d\le x} \psi_0(d)
\Big( g(d)|{\Cal O}_{\uL}(x/d,w)|^2-g(d+1) 
& |{\Cal O}_{\uL}(x/(d+1),w)|^2 \Big) \\
&- \psi_0(D) g(D+1)   |{\Cal O}_{\uL}(x/(D+1),w)|^2.
\tag{6.4}  \\
\endalign
$$
 Now we use that for $d>D$, $\psi_0(d) =O(d \exp(-(\log \log x)^2) =O(d/(\log X)^{A+2})$.  Hence 
 the second term in (6.4) is $\ll D^2 (\log X)^{-1} |{\Cal O}_{\uL}(x/(D+1),w)|^2 \ll x^2$ 
 upon using Proposition 5.4.
 The first term in (6.4) is 
 $$ 
\align
 \ll \sum_{D<d\le x} \frac{d}{(\log X)^{A+2}} &\Big( g ( d) \Big| |{\Cal O}_{\uL}(x/d,w)|^2 - 
 |{\Cal O}_{\uL}(x/(d+1),w)|^2 \Big| \\
&  + | g(d+1) -g(d)| |{\Cal O}_{\uL}(x/(d+1),w)|^2\Big).
 \\
 \endalign
 $$
 By Proposition 5.4 the second term above contributes $\ll x^2$, 
 while by Lemma 6.2 we have that the first term above is 
 $$ 
 \align
& \ll \sum_{D<d\le x} \frac{d^2}{(\log X)^A} \frac{x}{d}  \sum_{j=0}^{LR} 
w^j \sum_{x/((d+1)w^j)< n\le x/(dw^j)} |f(n)| 
 \\
 &\ll \frac{x}{ (\log X)^{A}} \sum_{j=0}^{LR} w^j \sum_{n\le x/w^j} |f(n)| \frac{x}{nw^j}\ll x^2, \\
 \endalign
 $$
 where the final estimate follows from Lemma 3.2.
   We conclude that the second term in the RHS of (6.3) is $\ll x^2$.

We  now turn to the first term in (6.3).  For any $D<d\le x$ and $d-1\le t\le d$ we have 
that $g(d)= g(t) + O(\log x)$, and by Lemma 6.2  that
$$
|{\Cal O}_{\uL} (x/d,w)|^2 =|{\Cal O}_{\uL}(x/t,w)|^2+ O\Big(\frac{x}{d} (\log X)^{\epsilon} 
\sum_{j=0}^{LR} w^j \sum_{x/(dw^j)<n\le x/((d-1)w^j)} |f(n)| \Big).
$$ 
Therefore, using also Proposition 5.4,
$$ 
\align
g(d) |{\Cal O}_{\uL}(x/d,w)|^2 &= \int_{d-1}^{d} g(t) |{\Cal O}_{\uL}(x/t,w)|^2 
dt  + O\Big(\frac{x^2}{d^2} (\log X)^{1+2\epsilon}\Big)
\\
& + O\Big(x(\log X)^{1+\epsilon} \sum_{j=0}^{LR} 
w^j \sum_{x/(dw^j) < n \le x/((d-1)w^j)} |f(n)| \Big).\\
\endalign
$$ 
Summing this over all $D<d\le x$ we get from the main term above the contribution 
$$ 
\int_D^x  t \log(ex/t) |{\Cal O}_{\uL}(x/t,w)|^2 dt \le x^2 \int_{1}^{x} |{\Cal O}_{\uL}(y,w)|^2 \log (ey) 
\frac{dy}{y^3}
$$ 
The error terms contribute 
$$
\align
&\ll \frac{x^2}{D} (\log X)^{1+2\epsilon} + x (\log X)^{1+\epsilon} 
\sum_{j=0}^{LR} w^j \sum_{n\le x/(Dw^j)} |f(n)|\\
& \ll x^2 + \frac{x^2}{D} (\log X)^{1+\epsilon} \sum_{j=0}^{LR} 
\sum_{n\le x/(Dw^j)} \frac{|f(n)|}{n} \ll x^2.\\ 
\endalign
$$
 The proof of the Proposition is complete.

\enddemo

We must now analyze the integral appearing in Proposition 6.3.    
To this end, we write 
$$ 
{\tilde S}(x) = \sum_{n\le x} f(n) \log n, 
$$ 
and define 
$$ 
{\tilde {\Cal  O}}_{\uell}(x,w) = \sum_{\uj \le \uell} (-1)^{j_1+\ldots+j_R} \binom{\uell}{\uj} w^{j_1(1+i\tau_1)+
\ldots+j_R(1+i\tau_R)} {\tilde S}(x/w^{j_1+\ldots+j_R}).
$$

\proclaim{Lemma 6.4}  We have 
$$
\align
\Big(\int_1^x  |{\Cal O}_{\uL}(t,w)|^2 \log (et) \frac{dt}{t^3} \Big)^{\frac 12} 
&\ll \Big(\int_{1}^{x} |{\tilde {\Cal O}}_{\uL}(t,w)|^2 \frac{dt}{t^3 \log (et)} \Big)^{\frac 12}
+  (\log X)^{7\epsilon/8}. 
\\
\endalign
$$
\endproclaim 
\demo{Proof}   We start as in the proof of Lemma 6.1.  Thus 
we may write 
$$ 
\align
&(\log t) {\Cal O}_{\uL}(t,w)= - \sum_{k=1}^{R} L w^{1+i\tau_k} (\log w) {\Cal O}_{\uL-\ue_k}(t/w,w)\\
&+\sum_{\uj \le \uL} (-1)^{j_1+\ldots+j_R} 
\binom{\uL}{\uj} \log (t/w^{j_1+\ldots+j_R}) S(t/w^{j_1+\ldots+j_R}) 
w^{j_1(1+i\tau_1)+\ldots+j_R(1+i\tau_R)}. 
\\
\endalign
$$ 
Since 
$$ 
(\log z)S(z) = {\tilde S}(z) + \sum_{n\le z} f(n) \log (z/n) = {\tilde S}(z) +\int_1^z \frac{S(y)}{y} dy,
$$
the second term above may be written as 
$$
{\tilde {\Cal O}}_{\uL} (t,w) + \int_1^t \frac{{\Cal O}_{\uL}(y,w)}{y} dy.
$$ 
Putting these remarks together, and using Proposition 5.4 we conclude that 
$$ 
(\log t) {\Cal O}_{\uL}(t,w) = {\tilde {\Cal O}}_{\uL}(t,w) + O(t (1+\log w)(\log X)^{2\epsilon/3}) 
= {\tilde {\Cal O}}_{\uL}(t,w) + O(t (\log X)^{5\epsilon/6}). 
$$ 
 The Lemma follows.  
 \enddemo


Putting together Proposition 6.3 and Lemma 6.4, we have 
that 
$$
|{\Cal O}_{\uL}(x,w)| \ll \frac{x}{\log x} (\log \log x)^{\frac 12} 
\Big( \int_{1}^{x} |{\tilde {\Cal O}}_{\uL}(t,w)|^2 \frac{dt}{t^3 \log (et)} \Big)^{\frac 12} 
 +  \frac{x}{\log x} (\log X)^{\epsilon} . 
 $$ 
Theorem 2.1 will now follow from the following Proposition.  

\proclaim{Proposition 6.5}  We have 
$$ 
\int_1^x  |{\tilde {\Cal O}}_{\uL}(t,w)|^2 \frac{dt}{t^3 \log (et)} 
\ll  (\log X)^{3\epsilon/2}. 
$$ 
\endproclaim
\demo{Proof}  We make the substitution $t=e^y$, obtaining
$$
\align
\int_1^x |{\tilde {\Cal O}}_{\uL}(t,w)|^2 \frac{dt}{t^3 \log (et)} 
&= \int_0^{\log x} |{\tilde {\Cal O}}_{\uL}(e^y,w) |^2e^{-2y} \frac{dy}{1+y}\\
&\ll \int_{1/\log X}^{\infty} e^{-2\alpha} \int_0^{\infty} |{\tilde {\Cal O}}_{\uL}(e^y,w)|^2 e^{- 2y(1+\alpha)} 
dy \ d\alpha. \tag{6.5}\\
\endalign
$$  
The idea now is to estimate the integral over $y$ in (6.5) using Plancherel's formula. 

Note that the Fourier transform of ${\tilde {\Cal  O}}_{\uL}(e^y,w) e^{-y(1+\alpha)}$ 
is 
$$ 
\align
\int_{-\infty}^{\infty} {\tilde {\Cal O}}_{\uL}(e^y,w) e^{-y(1+\alpha +it)} dy 
&= \sum_{\uj\le \uL} (-1)^{j_1+\ldots+j_R} \binom{\uL}{\uj} w^{j_1(1+i\tau_1)+\ldots+j_R(1+i\tau_R)} 
 \\
 &\times \int_{-\infty}^{\infty} 
\sum_{n\le e^y/w^{j_1+\ldots+j_R}} { f}(n) \log n  \ e^{-y(1+\alpha+it)} dy \\
&=\sum_{\uj \le \uL} (-1)^{j_1+\ldots+j_R} \binom{\uL}{\uj}
 w^{j_1(i\tau_1-\alpha-it)+\ldots+j_R(i\tau_R-\alpha-it)}\\
&\hskip .5 in \times 
\sum_{n=1}^{\infty} \frac{{ f}(n) \log n}{n^{1+\alpha+it}} \int_0^{\infty}  e^{-y(1+\alpha+it)} dy \\
&=- \frac{1}{(1+\alpha+it)}\prod_{k=1}^{R}
 (1-w^{-\alpha -it+i\tau_k})^{L} { F}^{\prime}(1+\alpha+it).\\
\endalign
$$
Therefore, by Plancherel's formula we have that  
$$ 
\align
\int_0^{\infty} |{\tilde {\Cal O}}_{\uL}(e^y,w)|^2 &e^{-y(2+2\alpha)} dy \\
&\ll  \int_{-\infty}^{\infty} |{F}^{\prime}(1+\alpha+it)|^2\prod_{k=1}^{R} 
\Big |1-w^{-\alpha-it+i\tau_k}\Big|^{2L} 
\frac{dt}{|1+\alpha+it|^{2}}.
\tag{6.6} \\
\endalign
$$

We split the integral in (6.6) into the regions when $|t|\le T/2$ and when $|t|>T/2$.  
For the latter region we use Lemma 3.2 which gives $|F^{\prime}(1+\alpha+it)| \ll (\log X)^{A+1}$ 
(note that   $\alpha\ge 1/\log X$ in (6.5)), so that this integral contributes 
$$
\ll (\log X)^{2A+2} \int_{|t|>T/2} \frac{dt}{|1+\alpha+it|^2} \ll 1.
$$
For the first region we use Lemmas 5.2 and 5.3 to obtain that 
$$
\align
|F^{\prime}(1+\alpha+it)| \prod_{k=1}^{R} |1-w^{-\alpha-it+i\tau_k}|^{L} 
&\ll \Big| \frac{F^{\prime}}{F}(1+\alpha+it) \Big|  |{\Cal F}_{\uL}(1+\alpha+it) | 
\\
&\ll \Big| \frac{F^{\prime}}{F}(1+\alpha+it) \Big| (\log X)^{\epsilon/2}.
\\
\endalign
$$ 
Therefore the integral over the first region contributes 
$$ 
\align
&\ll (\log X)^{\epsilon} \int_{|t|\le T/2} \Big|\frac{F^{\prime}}{F}(1+\alpha+it)\Big|^2 \frac{dt}{|1+\alpha+it|^2} 
\\
&\ll (\log X)^{\epsilon} \int_{-\infty}^{\infty} \Big|\frac{F^{\prime}}{F}(1+\alpha+it)\Big|^2 
\frac{dt}{|1+\alpha+it|^2}. \tag{6.7}
\\
\endalign
$$
Now the Fourier transform of the function $e^{-y(1+\alpha)} \sum_{n\le e^y} \Lam_f(n)$ 
is 
$$ 
\int_{-\infty}^{\infty} \sum_{n\le e^y} \Lam_f(n) e^{-y(1+\alpha+it)} dy 
= \sum_{n} \frac{\Lam_f(n)}{n^{1+\alpha+it}} \frac{1}{(1+\alpha+it)} 
= -\frac{F^{\prime}}{F}(1+\alpha+it) \frac{1}{(1+\alpha+it)},
$$ 
and so using Plancherel once again we obtain that the quantity in (6.7) is 
$$ 
\ll  (\log X)^{\epsilon} \int_{0}^{\infty} 
\Big| \sum_{n\le e^y} \Lam_f(n)\Big|^2 e^{-(2+2\alpha)y} dy. 
$$ 
We conclude that 
$$ 
\int_{0}^{\infty} |{\tilde {\Cal O}}_{\uL} (e^y,w)|^2 e^{-y(2+2\alpha)} dy \ll 1 + 
(\log X)^{\epsilon} \int_0^{\infty} \Big| \sum_{n\le e^y} \Lam_f(n)\Big|^2 e^{-(2+2\alpha)y}dy. 
\tag{6.8}
$$ 

Injecting the bound (6.8) in (6.5) we obtain that 
$$ 
\int_1^x |{\tilde {\Cal O}}_{\uL}(t,w)|^2 \frac{dt}{t^3 \log (et)} 
\ll 1 + (\log X)^{\epsilon} \int_{1/\log X}^{\infty} e^{-2\alpha} 
\int_0^{\infty} \Big| \sum_{n\le e^y} \Lam_f(n)\Big|^2 e^{-(2+2\alpha)y} dy \ d\alpha. 
$$
Expanding, we obtain that the double integrals above are 
$$ 
\align
&\ll \sum_{2\le n_1 \le n_2} |\Lam_f(n_1) \Lam_f(n_2)| \int_{1/\log X}^{\infty} \int_{\log n_2}^{\infty} e^{-(2+2\alpha)y-2\alpha} dy \ d\alpha \\
&\ll \sum_{2\le n_1 \le n_2} |\lam_f(n_1) \lam_f(n_2)| \frac{\Lam(n_1)\Lam(n_2)}{n_2^{2+2/\log X} \log n_2} \\ 
&\ll \sum_{2\le n_1 \le n_2}(|\lam_f(n_1)|^2 +|\lam_f(n_2)|^2)  \frac{\Lam(n_1) \Lam(n_2)}{n_2^{2+2/\log X}\log n_2} \\
&\ll \sum_{2\le n} \frac{|\lam_f(n)|^2\Lam(n)}{n^{1+2/\log X}\log n} \ll (\log \log X),\\
\endalign
$$ 
upon using the prime number theorem in the penultimate step, and (2.6)  for the last step. 
The Proposition follows, and with it Theorem 2.1. 

\enddemo

\enddocument